\documentclass[12pt,reqno]{amsart}

\usepackage{amssymb}
\usepackage{verbatim}
\usepackage{ifthen}

\renewcommand{\phi}{\varphi}

\newcommand{\Z}{{\mathbb Z}}

\newcommand{\cP}{\mathcal P}

\newcommand{\oA}{{\overline A}}
\newcommand{\oB}{{\overline B}}
\newcommand{\oC}{{\overline C}}
\newcommand{\oF}{{\overline F}}
\newcommand{\oS}{{\overline S}}

\newcommand{\oa}{{\overline a}}
\newcommand{\ob}{{\overline b}}
\newcommand{\of}{{\overline f}}
\newcommand{\oz}{{\overline z}}

\newcommand{\tA}{{\widetilde A}}
\newcommand{\tB}{{\widetilde B}}
\newcommand{\tC}{{\widetilde C}}

\newcommand{\pgh}[1][G/H]{\phi_{#1}}

\newcommand{\seq}{\subseteq}

\newcommand{\stm}{\setminus}
\newcommand{\est}{\varnothing}

\newcommand{\longc}{,\ldots,}
\newcommand{\longp}{+\dotsb+}

\newcommand{\division}[1]{\smallskip\noindent #1.\,}

\theoremstyle{plain}
\newtheorem{lemma}{Lemma}
\newtheorem{theorem}{Theorem}
\newtheorem{corollary}{Corollary}
\newtheorem{proposition}{Proposition}

\newtheorem{theirtheorem}{Theorem}

\theoremstyle{remark}
\newtheorem*{remark}{Remark}

\newcommand{\refl}[1]{~\ref{l:#1}}
\newcommand{\reft}[1]{~\ref{t:#1}}
\newcommand{\refc}[1]{~\ref{c:#1}}
\newcommand{\refp}[1]{~\ref{p:#1}}
\newcommand{\refs}[1]{~\ref{s:#1}}
\newcommand{\refb}[1]{~\cite{b:#1}}
\newcommand{\refe}[1]{~\eqref{e:#1}}

\title[Critical pairs and Kemperman's theorem]
  {Critical pairs in abelian groups \break and Kemperman's structure theorem}
\author{Vsevolod F. Lev}
\address{Department of Mathematics,
         The University of Haifa at Oranim,
         Tivon 36006,
         Israel}
\email{seva@math.haifa.ac.il}
\keywords{Critical pairs, small doubling, sumset, Kemperman's theorem, Kneser's theorem}
\subjclass[2000]{Primary: 11P70; secondary: 11B25, 11B13}

\begin{document}
\baselineskip 16pt

\begin{abstract}
A well-known result by Kemperman describes the structure of those pairs
$(A,B)$ of finite subsets of an abelian group satisfying $|A+B|\le|A|+|B|-1$.
We establish a description which is, in a sense, dual to Kemperman's, and as
an application sharpen several results due to Deshouillers, Hamidoune,
Hennecart, and Plagne.
\end{abstract}

\maketitle

\section{Overview of the paper}\label{s:overview}

The \emph{sumset} of two subsets $A$ and $B$ of an additively written group
is denoted by $A+B$ and defined as the set of all those group elements,
representable as a sum of an element of $A$ and an element of $B$:
  $$ A+B := \{ a+b\colon a\in A,\,b\in B \}. $$
In his remarkable paper \refb{k}, Kemperman classified completely all pairs
$(A,B)$ of finite subsets of an abelian group with the small sumset; more
precisely, those pairs satisfying
\begin{equation}\label{e:ApB}
  |A+B| \le |A|+|B|-1.
\end{equation}
This truly outstanding result, complementing a basic theorem by Kneser, may
have not received yet the recognition that it certainly deserves. One of the
reasons for this is that the structure of pairs, satisfying \refe{ApB}, is
rather complicated, and so is Kemperman's description reflecting this
structure. Indeed, \refb{k} is not an easy reading, and an essential part of
this paper constitutes an attempt to present Kemperman's theorem and
the mathematics around it in a possibly clear and transparent form, allowing
one to appreciate this highly non-trivial theorem and facilitating its
application.

To prepare the ground and explain where Kemperman's theorem stemmed from, we
start with the theorems of Kneser and Kemperman-Scherk; this is the subject
of the next section, where also the basic notion of period and some useful
notation are introduced.

In Section \refs{kempth} we first define \emph{elementary pairs}, which are
the ``building blocks'' of Kemperman's construction; once this is
accomplished, Kemperman's theorem is stated.

Having finished with the expository part, we present a number of new results
in Section \refs{summary}. Though the exact formulations are to be postponed
until the necessary notation is introduced and Kemperman's theorem (to which
our results are to be compared) is stated, we try to give here a brief
outline.

Our main result is Theorem \reft{dual}, which is in a sense dual to
Kemperman's theorem. In a very rough way, Kemperman's theorem can be compared
to the statement that if a pair $(a,b)$ of non-negative integers possesses
some property $\cP$, then there exist integers
 $d\ge 1,\,\oa,\ob\ge 0$, and $a_0,b_0\in[0,d-1]$ such that
$a=\oa d+a_0,\,b=\oa d+b_0$, and moreover, the pair $(a_0,b_0)$ is
``elementary'', while the pair $(\oa,\ob)$ satisfies $\cP$. Thus, we have a
closed-form description of the ``residual pair'' $(a_0,b_0)$ and a recursive
description of the ``root pair'' $(\oa,\ob)$. Theorem \reft{dual} can then be
compared to the statement that, under the same assumption, the integers
$d\ge1,\,\oa,\ob\ge 0$, and $a_0,b_0\in[0,d-1]$ with
 $a=\oa d+a_0,\,b=\oa d+b_0$ exist so that $(\oa,\ob)$ is elementary, while
$(a_0,b_0)$ satisfies $\cP$.

Another result established in this paper is Theorem \reft{me}, which can be
viewed as a light version of Kemperman's theorem. While the assumptions of
Theorem \reft{me} are slightly weaker than those of Kemperman's theorem, the
former provides somewhat less structure information than the later, giving
only a necessary (but not sufficient) condition for \refe{ApB} to hold. The
potential advantage of Theorem \reft{me} is its relative simplicity; at the
same time, it keeps much of the strength of Kemperman's theorem in the sense
that the deduction of each of these two results from another one is
reasonably short. Indeed, we do not have an independent proof of Theorem
\reft{me}; instead, we derive it from Kemperman's theorem, and then derive
Theorem \reft{dual} from Theorem \reft{me}.

As applications of Kemperman's theorem and Theorem \reft{dual}, we improve
several results by Deshouillers, Hamidoune, Hennecart, and Plagne; see
Section \refs{summary} for more information and discussion.

The proofs of all the results, presented in Section \refs{summary}, are
concentrated in Section \refs{proofs}. In the Appendix we give two
claims persistent to the definition of Kemperman's elementary pairs.

\section{Preliminaries: theorems of Kneser and Kemperman-Scherk}
        \label{s:KnKSth}

Let $G$ be an (additively written) abelian group. For finite non-empty
subsets $A,B\seq G$ with given cardinalities $|A|$ and $|B|$, how small the
sumset $A+B$ can be? It is not difficult to see that if $G$ is torsion-free,
then $|A+B|\ge|A|+|B|-1$ holds for any choice of $A$ and $B$. The well-known
Cauchy-Davenport theorem asserts that if $G$ is cyclic of prime order, then
this estimate remains true under the additional assumption that $A+B\neq G$.
Kneser's theorem is a deep and far-going extension of these results onto
arbitrary abelian groups; loosely speaking, it says that if $|A+B|<|A|+|B|-1$
holds, then some torsion is involved. To state the theorem precisely we have
to introduce some notation.

For a subset $S$ and an element $g$ of an abelian group $G$, we use the
abbreviations $S+g$ and $g+S$ for the sumset $S+\{g\}$. If $H$ is a subgroup
of $G$, then the canonical homomorphism of $G$ onto the quotient group $G/H$
is denoted by $\pgh$, and the full inverse image under $\pgh$ of a subset
$\oS\seq G/H$ is denoted by $\pgh^{-1}(\oS)$. The period of a subset $S\seq
G$ will be denoted by $\pi(S)$; recall, that this is the subgroup of $G$
defined by
  $$ \pi(S) := \{g\in G\colon S+g=S \}, $$
and that $S$ is called \emph{periodic} if $\pi(S)\neq\{0\}$ and
\emph{aperiodic} otherwise. Thus, $S$ is a union of cosets of $\pi(S)$ and
indeed, $\pi(S)$ lies above any subgroup of $G$ such that $S$ is a union of
its cosets. Observe also that if $H=\pi(S)$, then $\pgh(S)$ is an aperiodic
subset of $G/H$.

\begin{theirtheorem}[Kneser, \cite{b:kn1,b:kn2}; see also \refb{mann}]%
                                                           \label{t:kneser}
Let $A$ and $B$ be finite, non-empty subsets of an abelian group $G$,
satisfying \refe{ApB}. Then, letting $H := \pi(A+B)$, we have
  $$ |A+B| = |A+H|+|B+H|-|H|. $$
\end{theirtheorem}

This theorem will be referred to as ``Kneser's theorem''.
It shows that any pair $(A,B)$, satisfying \refe{ApB}, can be obtained by
``lifting'' a pair $(\oA,\oB)$ of subsets of a quotient group with
$|\oA+\oB|=|\oA|+|\oB|-1$. Indeed, if $A,B$, and $H$ are as in Kneser's
theorem, $\oA=\pgh(A)$, and $\oB=\pgh(B)$, then
$|A+H|=|\oA||H|,\,|B+H|=|\oB||H|$, and $|A+B|=|\oA+\oB||H|$, so that the
conclusion of the theorem takes the shape $|\oA+\oB|=|\oA|+|\oB|-1$.
Furthermore,
  $$ (|A+H|-|A|)+(|B+H|-|B|) = |H|+(|A+B|-|A|-|B|) < |H| $$
whence $A$ and $B$ are obtained from $A+H=\pgh^{-1}(\oA)$ and
$B+H=\pgh^{-1}(\oB)$ by removing less than $|H|$ elements totally.
Conversely, one verifies easily that if $H$ is a finite subgroup of an
abelian group $G$, the finite non-empty subsets $\oA,\oB\seq G/H$ satisfy
$|\oA+\oB|=|\oA|+|\oB|-1$, and $A\seq\pgh^{-1}(\oA),\,B\seq\pgh^{-1}(\oB)$
satisfy $(|\pgh^{-1}(\oA)|-|A|)+(|\pgh^{-1}(\oB)|-|B|)<|H|$, then \refe{ApB}
holds.

We need the following, almost immediate, corollary from Kneser's theorem.
\begin{corollary}\label{c:A/2+B}
Let $A$ and $B$ be finite, non-empty subsets of an abelian group, satisfying
$|A+B|<|A|/2+|B|$. Then $A$ is contained in a coset of $\pi(A+B)$.
\end{corollary}

\begin{proof}
Write $H:=\pi(A+B)$. Then by Kneser's theorem and the assumption,
  $$ |A+H|+|B+H|-|H| = |A+B| < |A+H|/2 + |B+H|, $$
and hence $|A+H|<2|H|$. Thus in fact $|A+H|=|H|$, as wanted.
\end{proof}

For finite subsets $A$ and $B$ of an abelian group and a group element $c$ we
write
  $$ \nu_c(A,B) := \# \{ (a,b)\in A\times B\colon c=a+b \}; $$
that is, $\nu_c(A,B)$ is the number of representations of $c$ as a sum of an
element of $A$ and an element of $B$. Observe, that $\nu_c(A,B)>0$ if and
only if $c\in A+B$. The smallest number of representations of an element of
$A+B$ will be denoted by $\mu(A,B)$:
  $$ \mu(A,B) := \min \{ \nu_c(A,B) \colon c\in A+B \}. $$

Another result frequently used in this paper is a beautiful theorem by
Kemperman and Scherk, relating the quantities $|A+B|$ and $\mu(A,B)$. We call
it ``the Kemperman-Scherk theorem''. This name is not standard; the reader
can check \hbox{\cite[Section~1]{b:l}} for explanations and historical
references.

\begin{theirtheorem}[Kemperman-Scherk]\label{t:ks}
If $A$ and $B$ are finite, non-empty subsets of an abelian group, then
  $$ |A+B| \ge |A| + |B| - \mu(A,B). $$
\end{theirtheorem}

Thus, the theorem of Kemperman-Scherk shows that if the sumset $A+B$ is
small, then any element $c\in A+B$ has many representations of the form
$c=a+b$ with $a\in A$ and $b\in B$.

Here is a short deduction of the theorem of Kemperman-Scherk from Kneser's
theorem.
\begin{proof}[Proof of Theorem \reft{ks}]
Let $H:=\pi(A+B)$, choose arbitrarily $c\in A+B$, and fix a representation
$c=a_0+b_0$ with $a_0\in A$ and $b_0\in B$. If \refe{ApB} fails, the
assertion is trivial; otherwise, by Kneser's theorem we have
\begin{align*}
  |(A-a_0)\cap H|+|(b_0-B)\cap H|
       &= |A\cap(a_0+H)|+|B\cap(b_0+H)| \\
       &= 2|H|-|(a_0+H)\stm A|-|(b_0+H)\stm B| \\
       &\ge 2|H|-|(A+H)\stm A|-|(B+H)\stm B| \\
       &= |H|+(|A|+|B|-|A+B|)
\end{align*}
and by the boxing principle, the sets $(A-a_0)\cap H$ and $(b_0-B)\cap H$
share at least $|A|+|B|-|A+B|$ common elements. However, each pair
 $(a,b)\in A\times B$ with $a-a_0=b_0-b$ yields the representation $c=a+b$,
whence $\nu_c(A,B)\ge |A|+|B|-|A+B|$.
\end{proof}

\section{Introduction: Kemperman's theorem}\label{s:kempth}

As we saw above, Kneser's theorem reduces the problem of classifying all
pairs $(A,B)$ of finite, non-empty subsets of an abelian group, satisfying
\refe{ApB}, to that of describing those pairs for which equality is attained
in \refe{ApB}. This latter problem was solved by Kemperman in \refb{k}. As
Kemperman has shown, one has to bring into consideration the following
additional condition:
\begin{equation}\label{e:kemcon}
  \text{either}\ \pi(A+B)=\{0\},\ \text{or}\ \mu(A,B)=1.
\end{equation}
Kneser's theorem shows that this condition is not too restrictive; more
precisely, if \refe{ApB} holds while Kemperman's condition \refe{kemcon}
fails, then there exist a finite non-zero subgroup $H\le G$ and subsets
$\oA,\oB\seq G/H$ such that $A$ and $B$ are ``dense'' in $\pgh^{-1}(\oA)$
and $\pgh^{-1}(\oB)$, respectively, and $(\oA,\oB)$ satisfies the analogs
of both \refe{ApB} and \refe{kemcon}. Furthermore, from theorems of
Kneser and Kemperman-Scherk it follows that if $(A,B)$ satisfies
\refe{ApB} and \refe{kemcon}, then in fact equality holds in \refe{ApB}.

Kemperman's description relies on the notion of an elementary pair, which
we introduce after short preparations.

By an arithmetic progression in an abelian group $G$ with difference
$d\in G$ we mean a set of the form $\{g+d,g+2d\longc g+kd\}$, where $g$
is an element of $G$ and $k$ is a positive integer, not exceeding the
order of $d$ in $G$. Thus, cosets of finite cyclic subgroups (and in
particular, one-element sets) are considered arithmetic progressions,
while the empty set is not.

For a subset $A$ of an abelian group we write $-A:=\{-a\colon a\in A\}$,
and given yet another subset $B$ and an element $g$ of the same group we
let $B-A:=B+(-A)$ and $g-A:=\{g\}-A$.

Following Kemperman, we say that a pair $(A,B)$ of finite subsets of an
abelian group $G$ is \emph{elementary} if at least one of the following
holds:
\begin{itemize}
\item[(I)] $\min\{|A|,|B|\}=1$;
\item[(II)] $A$ and $B$ are arithmetic progressions sharing a common
  difference $d\in G$, the order of which in $G$ is at least $|A|+|B|-1$;
\item[(III)] $A=g_1+(H_1\cup\{0\})$ and $B=g_2-(H_2\cup\{0\})$, where
  $g_1,g_2\in G$ and $H_1,H_2$ are non-empty subsets of a subgroup
  $H\le G$ such that $H=H_1\cup H_2\cup\{0\}$ is a partition of $H$.
  Moreover, $c:=g_1+g_2$ is the only element of $A+B$ with $\nu_c(A,B)=1$;
\item[(IV)] $A=g_1+H_1$ and $B=g_2-H_2$, where $g_1,g_2\in G$ and $H_1,H_2$
  are non-empty, aperiodic subsets of a subgroup $H\le G$ such that
  $H=H_1\cup H_2$ is a partition of $H$. Moreover, $\mu(A,B)\ge 2$.
\end{itemize}

We note that if $(A,B)$ is an elementary pair of subsets, then $A$ and $B$
are non-empty; also, the subgroup $H$ in (III) and (IV) is finite and
satisfies $|H|\ge 3$. Another immediate, yet important observation is that if
$A$ and $B$ are subsets of a subgroup $F\le G$, then $(A,B)$ is an elementary
pair in $G$ of some type if and only if it is an elementary pair in $F$ of
the same type.

Notice, that for type (III) pairs we have $|A|+|B|=|H|+1$ whence
$A+B=g_1+g_2+H$ by the boxing principle, while for type (IV) pairs we have
$|A|+|B|=|H|$ and $A+B=g_1+g_2+(H\stm\{0\})$. (To prove the last equality fix
an element $h\in H\stm\{0\}$ and notice that since $h\notin\pi(H_2)$, there
exist $h_2\in H_2$ and $h_1\in H_1$ with $h+h_2=h_1$, whence
$g_1+g_2+h=(g_1+h_1)+(g_2-h_2)\in A+B$.) One now verifies easily that both
\refe{ApB} and \refe{kemcon} hold true for elementary pairs $(A,B)$ of any
type.

It can be shown that the first sentence in (III) describes precisely those
pairs $(A,B)$, satisfying both \refe{ApB} and \refe{kemcon} and such that
$A+B$ is an $H$-coset; the second sentence excludes the situation where
$(A,B)$ can be decomposed into ``more elementary'' pairs. Similarly, the
first sentence in (IV) describes precisely those pairs $(A,B)$ satisfying
both \refe{ApB} and \refe{kemcon} and such that $A+B$ is an $H$-coset with
one element removed, while the second one excludes the situation where
$(A,B)$ can be further decomposed. We make these statements precise and prove
them in the Appendix.

Now we can state the result which throughout the rest of the paper is
referred to as ``Kemperman's theorem''.
\begin{theirtheorem}[Kemperman, \protect{\cite[Theorem~5.1]{b:k}}]%
                                                              \label{t:KMST}
Let $A$ and $B$ be finite, non-empty subsets of a non-trivial abelian group
$G$. A necessary and sufficient condition for both \refe{ApB} and
\refe{kemcon} to hold simultaneously is that there exist non-empty subsets
$A_0\seq A,\,B_0\seq B$ and a non-zero subgroup $H\le G$ such that, writing
$\oA:=\pgh(A)$ and $\oB:=\pgh(B)$, we have
\begin{itemize}
\item[(i)]   each of $A_0$ and $B_0$ is contained in an $H$-coset, and the
  pair $(A_0,B_0)$ is elementary;
\item[(ii)]  each of $A\stm A_0$ and $B\stm B_0$ is a (possibly, empty) union
  of $H$-cosets;
\item[(iii)] $|\oA+\oB|=|\oA|+|\oB|-1$;
\item[(iv)]  $A_0+B_0+H$, considered as an element of $G/H$, has a unique
  representation as a sum of an element of $\oA$ and an element of $\oB$.
\end{itemize}
\end{theirtheorem}

We notice that one can choose $H=G$ in Kemperman's theorem if and only if
$(A,B)$ is an elementary pair. Observe also that, as it follows from (ii),
the set $(A+B)\stm(A_0+B_0)$ is a (possibly, empty) union of $H$-cosets.
Consequently, if $H$ is finite then conditions (i) and (ii) show that
$|A+H|-|A|\le |H|-1,\,|B+H|-|B|\le |H|-1$, and $|A+B+H|-|A+B|\le |H|-1$
(cf.~Theorem \reft{me} below). Note, that for a finite set $C\seq G$ and a
finite subgroup $H\le G$ the condition $|C+H|-|C|\le |H|-1$ (frequently
emerging in what follows, and sometimes written as $|C+H|\le |C|+|H|-1$)
means that $C$ is a ``dense'' subset of a union of $H$-cosets.

\section{Summary of results}\label{s:summary}

In this section we present and discuss our results; the proofs are postponed
until Section \refs{proofs}.

While Kemperman's theorem describes the structure of the ``residual'' pair
$(A_0,B_0)$, in the following theorem (qualified above as a light version of
Kemperman's theorem) we make the emphasis on the structure of the ``root''
pair $(\oA,\oB)$. Notice, that Kemperman's condition \refe{kemcon} is
replaced by a weaker assumption; on the other hand, we obtain a less precise
information about the distribution of elements of $A$ and $B$ in $H$-cosets
(cf. Theorem \reft{dual} below).
\begin{theorem}\label{t:me}
Let $A$ and $B$ be finite, non-empty subsets of a non-trivial abelian group
$G$, satisfying \refe{ApB}. Suppose that either $A+B\neq G$ or $\mu(A,B)=1$
(or both) hold true. Then there exists a finite, proper subgroup $H<G$ such
that, writing $\oA=\pgh(A)$ and $\oB=\pgh(B)$, we have
\begin{itemize}
\item[(i)]  $|C+H|-|C|\le|H|-1$ with $C$ substituted by any of the sets $A,B$,
  and $A+B$;
\item[(ii)] $(\oA,\oB)$ is an elementary pair in $G/H$.
\end{itemize}
\end{theorem}

Theorem \reft{me} shows that, under some mild technical restriction, any
pair of finite subsets $(A,B)$ of an abelian group $G$, satisfying
\refe{ApB}, can be obtained by lifting an elementary pair $(\oA,\oB)$ of
subsets of a quotient group $G/H$ and removing ``a small number'' of elements
from the lifted subsets. Concerning the technical restriction just
mentioned, we observe that the condition ``either $A+B\neq G$ or $\mu(A,B)=1$
(or both) hold true'' means that there is a group element with at most one
representation as $a+b$ with $a\in A$ and $b\in B$.

Our next theorem provides an alternative classification of pairs $(A,B)$,
satisfying \refe{ApB}. The reader can now appreciate the statement made in
Section \refs{overview} on duality between this theorem and Kemperman's one.
\begin{theorem}\label{t:dual}
Let $A$ and $B$ be finite, non-empty subsets of an abelian group $G$. A
necessary and sufficient condition for both \refe{ApB} and \refe{kemcon} to
hold simultaneously is that either the pair $(A,B)$ is elementary, or there
exist non-empty subsets $A_0\seq A,\ B_0\seq B$ and a finite, non-zero,
proper subgroup $H<G$ such that, writing $\oA:=\pgh(A)$ and $\oB:=\pgh(B)$,
we have
\begin{itemize}
\item[(i)]   each of $A_0,B_0$ is contained in an $H$-coset,
  $|A_0+B_0|=|A_0|+|B_0|-1$, and $(A_0,B_0)$ satisfies Kemperman's condition
  \refe{kemcon};
\item[(ii)]  each of $A\stm A_0$ and $B\stm B_0$ is a (possibly, empty) union
  of $H$-cosets;
\item[(iii)] the pair $(\oA,\oB)$ is elementary, and $A_0+B_0+H$, considered
  as an element of $G/H$, has a unique representation as a sum of an
  element of $\oA$ and an element of $\oB$.
\end{itemize}
\end{theorem}

We point it out once again that Theorem \reft{dual} will be derived from
Theorem \reft{me} which, in turn, is deduced from Kemperman's theorem.

Turning to applications, we first show how Theorem \reft{me} can be used to
sharpen one of the central results of \refb{hp2}.

We say that a non-empty, finite subset $A$ of an abelian group $G$ is a
\emph{thick component} of $G$ (the term used in \refb{hp2} is ``Vosper
subset'') if for any finite subset $B\seq G$ with $|B|\ge 2$ either
$|A+B|\ge|A|+|B|$, or $|A|+|B|\ge|G|-1$ (or both) hold. Evidently, if $|G|>4$
and $A$ is a thick component of $G$, then $|A|\ge 2$.
\begin{theirtheorem}[\protect{\cite[Theorem~2.1]{b:hp2}}]
Let $A$ be a generating subset of a finite abelian group $G$. Suppose that
$0\in A$ and $|A|\le|G|/2$. Then there exists a subgroup $H<G$ satisfying
$|A+H|-|A|\le |H|-1$ and $A+H\neq G$, and such that $\pgh(A)$ is either an
arithmetic progression or a thick component (of $G/H$).
\end{theirtheorem}

Our version is as follows.
\begin{theorem}\label{t:thick}
Let $A$ be a finite, non-empty subset of an abelian group $G$. Suppose that
$|A|\le|G|/2$. Then there exists a finite subgroup $H<G$ satisfying
$|A+H|-|A|\le |H|-1$ and $A+H\neq G$, and such that $\pgh(A)$ is either an
arithmetic progression or a thick component (of $G/H$).
\end{theorem}

Next, we show how the original Kemperman's theorem leads to an improvement of
a result of Hamidoune \cite[Theorem~6.6]{b:h} (see \cite[Lemma~2]{b:hp1} and
\cite[Theorem~A and Lemma~3.2]{b:hp2} for the corrected statement).
\begin{theorem}\label{t:hampla}
Let $A$ and $B$ be finite subsets of an abelian group $G$, satisfying
\refe{ApB} and such that $A$ is non-empty, $|B|\ge 2$, and $|A+B|\le
|G|-2$. Then at least one of the following holds:
\begin{itemize}
\item[(i)]  there exists a finite, non-zero subgroup $H<G$ such that
  $|A+H|-|A|\le |H|-1$ and $A+H\neq G$;
\item[(ii)] there exist a non-zero subgroup $H\le G$ of finite index
  $n:=[G:H]$, elements $g_1\longc g_{n-1}\in G$, and an arithmetic
  progression $A_0\seq A$, such that $A$ is the
  disjoint union $A = A_0\cup(g_1+H)\cup\dots\cup(g_{n-1}+H)$.
\end{itemize}
\end{theorem}

Observe that if the subgroup $H$ in (ii) satisfies $H=G$, then $A$ is an
arithmetic progression. Also, if $H$ is infinite then $H=G$ necessarily holds
since $A$ is finite.

\begin{remark}
As compared with Theorem \reft{hampla}, \cite[Theorem~A]{b:hp2} imposes
the extra restrictions that $G$ is finite, $A$ contains zero and
generates $G$, and $|A|\le|G|/2$.
\end{remark}

Finally, we apply Theorem \reft{dual} to the situation where $G$ is an
elementary abelian $2$-group, and the set summands $A$ and $B$ are identical.
That is, we are interested in the structure of those subsets
$A\seq(\Z/2\Z)^r$, with a positive integer $r$, satisfying $|2A|<2|A|$. (Here
and below, $2A$ is used as an abbreviation for $A+A$.) We show that Theorem
\reft{dual} allows one to give a complete, closed-form (non-recursive)
classification of such subsets. We start with two theorems, describing the
present state of the art; the former of them establishes the structure of the
sumset $2A$ (for a set $A\seq(\Z/2\Z)^r$ with $|2A|<2|A|$), the latter
presents some structure information about the set $A$ itself.
\begin{theirtheorem}[\protect{\cite[Theorem 1]{b:hnp}}]
Let $r\ge 1$ be an integer and suppose that $A$ is a non-empty subset of the
elementary abelian $2$-group $(\Z/2\Z)^r$, satisfying $|2A|=c|A|$ with $c<2$.
Then there exist a subgroup $H\le(\Z/2\Z)^r$ and an element $g\in(\Z/2\Z)^r$
such that $A\seq g+H$, and moreover,
\begin{itemize}
\item[(i)]  if $1\le c<7/4$ then $2A=H$;
\item[(ii)] if $7/4\le c<2$ then either $2A=H$, or there exist a subgroup
  $F<H$ with $|F|\le |H|/8$ and an element $h_0\in H$ such that
  $2A=H\stm(h_0+F)$.
\end{itemize}
\end{theirtheorem}

The restriction of \cite[Theorem~2]{b:dhp} onto the range $c\in[1,2)$ is
\begin{theirtheorem}[cf. \protect{\cite[Theorem~2]{b:dhp}}]
Let $r\ge 1$ be an integer and suppose that $A$ is a non-empty subset of the
elementary abelian $2$-group $(\Z/2\Z)^r$, satisfying $|2A|=c|A|$ with $c<2$.
Then there exists a subgroup $H\le(\Z/2\Z)^r$ such that $A\seq H$ and $|A|\ge
u(c)|H|$, where $u(c)=(3c-1-c^2)/(2c-1)$.
\end{theirtheorem}

We are now in a position to state the necessary and sufficient condition for
$|2A|<2|A|$ (with $A\seq(\Z/2\Z)^r$) to hold. Suppose that $H\le(\Z/2\Z)^r$
is a non-zero subgroup and let $h_0\in H$ be a non-zero element of $H$. We
say that $S\seq H$ is an $h_0$-antisymmetric subset of $H$, if $H$ is the
disjoint union $H=S\cup(h_0+S)$; in other words, $S$ contains exactly one
element from each pair $(h,h+h_0)$, for any $h\in H$. Evidently, in this case
we have $|S|=|H|/2$.

\begin{theorem}\label{t:struct}
Let $r\ge 1$ be an integer. If a subset $A\seq(\Z/2\Z)^r$ satisfies
$|2A|<2|A|$, then one of the following holds:
\begin{itemize}
\item[(i)]  there exists a subgroup $H\le(\Z/2\Z)^r$ such that $A$ is contained in
  an $H$-coset and $|A|>|H|/2$;
\item[(ii)] there exist two subgroups $F,H\le(\Z/2\Z)^r$, satisfying $|F|\ge 8$
and $F\cap H=\{0\}$, and an aperiodic antisymmetric subset $S\seq F$, such
that $A$ is obtained from a shift of the set $S+H$ by removing less than
$|H|/2$ of its elements. In this case $2A$ is the sum $F\oplus H$ with one
$H$-coset removed, so that $|2A|=(|F|-1)|H|$.
\end{itemize}

Conversely, let $F,H\le(\Z/2\Z)^r$ be subgroups satisfying $F\neq\{0\}$ and
$F\cap H=\{0\}$, and let $S$ be an antisymmetric subset of $F$. If $A$ is
obtained from a shift of the set $S+H$ by removing less than $|H|/2$ of its
elements, then $|2A|=(|F|-|\pi(S)|)|H|<2|A|$.
\end{theorem}

\section{Proofs}\label{s:proofs}

\begin{proof}[Proof of Theorem \reft{me}]
Suppose that $A$ and $B$ are finite, non-empty subsets of a non-trivial
abelian group $G$, satisfying \refe{ApB} and such that either $A+B\neq G$
or $\mu(A,B)=1$ (or both) hold true. We want to show that there exists a
finite proper subgroup $H<G$ with the properties that
\begin{itemize}
\item[(i)]  $|C+H|-|C|\le|H|-1$ for $C=A,B,A+B$;
\item[(ii)] $(\oA,\oB)$ is an elementary pair in $G/H$, where $\oA=\pgh(A)$
  and $\oB=\pgh(B)$.
\end{itemize}

If there exists a finite proper subgroup $H<G$ such that
$|A+B+H|-|A+B|\le|H|-1$ and $|H|\ge|A+B|$, then $|A+B+H|<2|H|$;
consequently, $A+B$, and therefore each of the sets $A$ and $B$, is
contained in a coset of $H$, from which the assertion follows.

Suppose now that for each finite proper subgroup $H<G$ with
$|A+B+H|-|A+B|\le|H|-1$ we have $|H|<|A+B|$. Observing that the zero subgroup
$H=\{0\}$ satisfies the conditions
\begin{itemize}
\item[(a)] $|C+H|-|C|\le|H|-1$ for $C=A,B,A+B$;
\item[(b)] either $\oA+\oB\neq G/H$ or $\mu(\oA,\oB)=1$ (or both) hold true;
\item[(c)] $|\oA+\oB|\le |\oA|+|\oB|-1$,
\end{itemize}
we can then find a finite proper subgroup $H<G$ which also satisfies these
conditions and which is maximal by inclusion among all finite proper
subgroups with this property. (Here and below in the proof $\oA$ and $\oB$
are defined by $\oA=\pgh(A),\,\oB=\pgh(B)$.) We show that the pair
$(\oA,\oB)$ is elementary.

Assume first that $(\oA,\oB)$ satisfies Kemperman's condition
\refe{kemcon}. Let $\oF\le G/H$ be the non-zero subgroup, the existence
of which is guaranteed by Kemperman's theorem as applied to the subsets
 $\oA,\oB\seq G/H$, and let $F:=\pgh^{-1}(\oF)$, so that $H\lneqq F$.
If $(\oA,\oB)$ is not elementary, then $\oF$ is finite and $\oF\lneqq G/H$,
whence $F$ is finite, too, and $F\lneqq G$. Writing $\tC=\pgh[G/F](C)$ for a
subset $C\seq G$, we show that $F$ satisfies the conditions
\begin{itemize}
\item[(a$'$)] $|C+F|-|C|\le|F|-1$ for $C=A,B,A+B$;
\item[(b$'$)] either $\tA+\tB\neq G/F$ or $\mu(\tA,\tB)=1$ (or both) hold true;
\item[(c$'$)] $|\tA+\tB|\le |\tA|+|\tB|-1$,
\end{itemize}
contradicting maximality of $H$.

Indeed, (c$'$) follows immediately from assertion (iii) of Kemperman's
theorem in view of the natural isomorphism $(G/H)/\oF\cong G/F$. Similarly,
(b$'$) follows from assertion (iv) of Kemperman's theorem, which implies that
$\mu(\tA,\tB)=1$. Finally, to prove (a$'$) we observe that for $C$ being any
of $A,B,A+B$ we have
  $$ (|\tC|-1)|\oF| = |\oC+\oF|-|\oF| \le |\oC|-1 $$
(see a remark following the statement of Kemperman's theorem). Multiplication
by $|H|$ and condition (a) yield
  $$ (|\tC|-1)|F| \le |C+H|-|H| \le |C|-1, $$
whence
  $$ |C+F|-|F| \le |C|-1, $$
as wanted.

Now assume that the pair $(\oA,\oB)$ does not satisfy Kemperman's condition
\refe{kemcon}. In this case $\oA+\oB$ is periodic and we define $\oF\le G/H$
to be the period: $\oF:=\pi(\oA+\oB)$; also, let $F:=\pgh^{-1}(\oF)$. By (b),
we have $\oF\neq G/H$, whence $H\lneqq F\lneqq G$. As above, for $C\seq G$ we
write $\tC=\pgh[G/F](C)$ and show that $F$ satisfies (a$'$)--(c$'$),
contradicting maximality of $H$.

By Kneser's theorem we have
  $$ |\oA+\oB|=|\oA+\oF|+|\oB+\oF|-|\oF|; $$
equivalently, $|\tA+\tB|=|\tA|+|\tB|-1$, establishing (c$'$). Next,
$\tA+\tB=G/F$ would imply $\oF=G/H$ which is wrong; this proves (b$'$).
Finally,
  $$ |A+B+F|=|A+B+H| \le |A+B|+|H|-1 < |A+B|+|F|-1 $$
(as $\oA+\oB+\oF=\oA+\oB$) and furthermore,
\begin{gather*}
  |\oA+\oF|+|\oB+\oF|-|\oF| = |\oA+\oB| \le |\oA|+|\oB|-1, \\
  |\oA+\oF| \le |\oA|+|\oF|-1, \\
  |A+F| \le |A+H|+|F|-|H| \le |A|+|F|-1;\\
\intertext{similarly,}
  |B+F| \le |B|+|F|-1.
\end{gather*}
This proves (a$'$).
\end{proof}

\begin{proof}[Proof of Theorem \reft{dual}]
Leaving the simple sufficiency part as an exercise to the interested reader,
we prove necessity: assuming that the pair $(A,B)$ satisfies \refe{ApB} and
\refe{kemcon} and is not elementary (and in particular $\min\{|A|,|B|\}>1$),
we establish the existence of $A_0,B_0$, and $H$ possessing properties
(i)--(iii). We divide the proof into two parts, the first of which will be
further subdivided.

\division{1} First, consider the situation where at least one of the sets $A$
and $B$ is contained in a coset of a finite proper subgroup of $G$. Suppose,
for instance, that $A$ is a subset of a coset of a finite proper subgroup
$H<G$, and let, moreover, $H$ be a minimal (by inclusion) subgroup with this
property. If $B$ is also contained in an $H$-coset, the assertion is
immediate; otherwise, we represent $B$ as a union
  $$ B=B_0\cup B_1\cup\dots\cup B_{n-1} $$
where $n\ge 2$, each subset $B_i$ is contained in a coset of $H$, and these
$n$ cosets are pairwise disjoint. Notice, that then $A+B$ is the disjoint
union
  $$ A+B=(A+B_0)\cup(A+B_1)\cup\dots\cup(A+B_{n-1}). $$

We recall that $(A,B)$ satisfies either $\pi(A+B)=\{0\}$, or $\mu(A,B)=1$.

\division{1.1} If $\mu(A,B)=1$ then, renumbering the sets $B_i$ if necessary,
we can assume that $\mu(A,B_0)=1$, whence $|A+B_0|\ge|A|+|B_0|-1$ by the
theorem of Kemperman-Scherk. Combining this with the trivial estimate
$|A+B_i|\ge|B_i|\ (i\in[1,n-1])$ we get
\begin{multline*}
  |A|+|B_0|+|B_1|\longp|B_{n-1}|-1 \ge |A+B| \\
                                      \ge (|A|+|B_0|-1)+|B_1|\longp|B_{n-1}|.
\end{multline*}
It follows that $|A+B_0|=|A|+|B_0|-1$ and $|A+B_i|=|B_i|$ for all indices
$i\in[1,n-1]$. The latter is only possible if $A$ is contained in a coset of
$\pi(B_i)$. By minimality of $H$, we have then $\pi(B_i)=H$; that is, $B_i$
is an $H$-coset for any $i\in[1,n-1]$. This establishes the required
structure (with $A_0=A$).

\division{1.2} If $\mu(A,B)>1$, then $\pi(A+B)=\{0\}$. In this case we first
observe that, by \refe{ApB}, there is at most one index $i\in[0,n-1]$ such
that $|A+B_i|\ge|B_i|+|A|/2$; without loss of generality, we can assume that
$|A+B_i|<|B_i|+|A|/2$ for $i\in[1,n-1]$. Corollary \refc{A/2+B} shows that
$A$ is contained in a coset of $\pi(A+B_i)$ for $i\in[1,n-1]$, and by
minimality of $H$, all $A+B_i\ (i\in[1,n-1])$ are $H$-cosets. Now, since
$A+B$ is aperiodic, so is $A+B_0$, whence $|A+B_0|\ge |A|+|B_0|-1$ by
Kneser's theorem. Finally,
\begin{multline*}
  |A|+|B_0|+|B_1|\longp|B_{n-1}|-1 \ge |A+B| \\
                               = |A+B_0|+(n-1)|H| \ge |A|+|B_0|+(n-1)|H|-1,
\end{multline*}
which implies $|A+B_0|=|A|+|B_0|-1$ and $|B_i|=|H|$ for $i\in[1,n-1]$.

\division{2} For the rest of the proof we assume that neither $A$ nor $B$ is
contained in a coset of a finite proper subgroup of $G$. We find then a
subgroup $H$ as in Theorem \reft{me}, so that $H$ is finite, $\{0\}\lneqq
H\lneqq G$ ($H$ is non-zero as the pair $(A,B)$ is not elementary), and
$(\oA,\oB)$ is elementary of one of the types (II)--(IV); here and below in
the proof $\oA$ and $\oB$ are defined by $\oA:=\pgh(A)$ and $\oB:=\pgh(B)$,
as in Theorem \reft{me}. Write
  $$ A = A_0\cup A_1\cup\dots\cup A_{m-1} \quad\text{and}\quad
                                  B = B_0\cup B_1\cup\dots\cup B_{n-1}, $$
where $m,n\ge 2$, each of the sets $A_i$ and $B_j$ is contained in an
$H$-coset, the cosets $A_i+H$ are pairwise disjoint, and so are the cosets
$B_j+H$. Renumbering the sets $A_i$ and $B_j$ we assume that if $\mu(A,B)=1$,
then there exists $c\in A+B$ with
\begin{equation}\label{e:incase1}
  \nu_c(A_0,B_0)=\nu_c(A,B)=1,
\end{equation}
and if $\mu(A,B)>1$ and $\pi(A+B)=\{0\}$ then
\begin{equation}\label{e:incase2}
  A_0+B_0+H\nsubseteq A+B.
\end{equation}
Furthermore, interchanging $A$ and $B$, if necessary, we can assume that
\begin{equation}\label{e:A0leB0}
  |A_0|\le|B_0|.
\end{equation}
We notice that for each $i\in[1,m-1]$ we have
\begin{equation}\label{e:A0Ai}
  |A_0|+|A_i|\ge 2|H|-(|A+H|-|A|) \ge |H|+1,
\end{equation}
and similarly,
\begin{equation}\label{e:B0Bi}
  |B_0|+|B_j|\ge |H|+1
\end{equation}
for each $j\in[1,n-1]$.

If $\mu(A,B)=1$ then $\mu(A_0,B_0)=1$ by \refe{incase1} and consequently
\begin{equation}\label{e:A0pB0issmall}
  |A_0|+|B_0|\le |H|+1
\end{equation}
by the boxing principle. If, on the other hand, $\mu(A,B)>1$ and
$\pi(A+B)=\{0\}$, then \refe{A0pB0issmall} (and in fact, the stronger
estimate $|A_0|+|B_0|\le |H|$) follows from the boxing principle and
\refe{incase2}. We see that \refe{A0pB0issmall} holds in any case, and
comparing it with \refe{A0Ai} and \refe{B0Bi} we conclude that
\begin{equation}\label{e:AiBjH}
  |A_i|+|B_j|\ge |H|+1
\end{equation}
for each $i\in[1,m-1]$ and $j\in[1,n-1]$. Also,
\begin{equation}\label{e:loc201}
  |A_i|+|B_0|\ge |A_i|+|A_0|\ge |H|+1
\end{equation}
for each $i\in[1,m-1]$ by \refe{A0leB0} and \refe{A0Ai}. From \refe{AiBjH}
and \refe{loc201} we conclude now that $A_i+B_j$ is an $H$-coset for any
$i\in[1,m-1]$ and $j\in[0,n-1]$.

It follows that $A_0+B_0+H\in G/H$ has a unique representation as a sum of an
element of $\oA$ and an element of $\oB$: for, assuming $A_i+B_j=A_0+B_0+H$
with $i\in[1,m-1]$ and $j\in[1,n-1]$, we obtain a contradiction with both
\refe{incase1} and \refe{incase2}. Thus, $(\oA,\oB)$ is type (II) or (III),
and in either case, each element of $\oA+\oB$, except $A_0+B_0+H$, has a
representation as $(A_i+H)+(B_j+H)$ with $i\in[1,m-1]$ and $j\in[0,n-1]$;
accordingly, $A+B$ is a disjoint union of $|\oA+\oB|-1=m+n-2$ full $H$-cosets
and the ``residual'' set $A_0+B_0$. The number of elements of this residual
set can now be easily estimated: if $\mu(A_0,B_0)=1$ then
\begin{equation}\label{e:A0B0islarge}
  |A_0+B_0|\ge |A_0| + |B_0| - 1
\end{equation}
by the theorem of Kemperman-Scherk, and if $A+B$ is aperiodic then so is
$A_0+B_0$, and hence \refe{A0B0islarge} holds true by Kneser's theorem. To
complete the proof we notice that
\begin{multline*}
  |A_0|+|A_1|\longp|A_{m-1}|+|B_0|+|B_1|\longp|B_{n-1}|-1 \\
  \ge |A+B| = (m+n-2)|H|+|A_0+B_0| \\
  \ge (m+n-2)|H|+|A_0|+|B_0|-1.
\end{multline*}
This shows that $|A_i|=|B_j|=|H|$ for all $i\in[1,m-1],\,j\in[1,n-1]$, and
that equality holds in \refe{A0B0islarge}.
\end{proof}

\begin{proof}[Proof of Theorem \reft{thick}]
If $A$ is contained in a coset of a finite proper subgroup of $G$ then we can
take this subgroup for $H$ and the result follows as $\pgh(A)$ is then a
one-element set, hence an arithmetic progression. For the rest of the proof
we assume that $A$ is \emph{not} contained in a coset of a finite proper
subgroup.

Suppose that $H<G$ is a finite subgroup such that
\begin{equation}\label{e:AHast}
  |A+H| - |A| \le |H| - 1 \quad\text{and}\quad A+H \neq G.
\end{equation}
Since $A$ is not contained in an $H$-coset, we derive from \refe{AHast} that
$2|H|\le |A+H|<|A|+|H|$, whence $|H|<|A|$. As $H=\{0\}$ satisfies
\refe{AHast}, there exists a finite subgroup $H<G$ satisfying \refe{AHast}
and maximal in the sense that it is not properly contained in any other
finite subgroup satisfying \refe{AHast}. We show that $\oA:=\pgh(A)$ is
either an arithmetic progression or a thick component of $G/H$.

Indeed, assume that $|\oA|\ge 2$ and that $\oA$ is \emph{not} a thick
component of $G/H$, so that $|\oA+\oB|<\min\{|G/H|-1, |\oA|+|\oB|\}$ for some
finite subset $\oB\seq G/H$ with $|\oB|\ge 2$. By Theorem \reft{dual}, either
$(\oA,\oB)$ is an elementary pair, or there exist a finite non-zero proper
subgroup $\oF<G/H$ and subsets $\oA_0\seq\oA$ and $\oB_0\seq\oB$ such that
both $\oA_0$ and $\oB_0$ are contained in $\oF$-cosets, and both
$\oA\stm\oA_0$ and $\oB\stm\oB_0$ are unions of $\oF$-cosets. In the former
case ($(\oA,\oB)$ is an elementary pair) from
 $\min\{|\oA|,|\oB|\}\ge 2,\ |\oA+\oB|\le|G/H|-2$, and the fact that $\oA$ is
not contained in a finite proper coset of $G/H$ it follows that $\oA$ is an
arithmetic progression. In the latter case, writing $F:=\pgh^{-1}(\oF)$, we
get
  $$ |A+F| = |\oA+\oF||H| \le (|\oA|+|\oF|-1)|H|
                                     = |A+H| + |F| - |H| \le |A| + |F| -1, $$
whence $A+F=G$ by maximality of $H$. Thus $\oA+\oF=G/H$ and therefore
\begin{gather*}
  |\oA|=|G/H|-|\oF|+|\oA_0|\ge|G/H|-|\oF|+1, \\
  |A|>|A+H|-|H|=(|\oA|-1)|H|\ge(|G/H|-|\oF|)|H|=|G|-|F|.
\end{gather*}
Now $|A|\le|G|/2$ gives $|F|>|G|/2$ and hence $F=G$, contradicting $\oF\lneqq
G/H$.
\end{proof}

\begin{proof}[Proof of Theorem \reft{hampla}]
Clearly, the period $P:=\pi(A+B)$ is finite and satisfies $A+P\neq G$ in view
of $A+B\neq G$. Furthermore, $|A+P|-|A|\le|P|-1$ by Kneser's theorem, and if
$P$ is non-zero then the assertion follows. Suppose now that $P=\{0\}$; that
is, $A+B$ is aperiodic. In this case the pair $(A,B)$ satisfies \refe{kemcon}
and we find subsets $A_0\seq A$ and $B_0\seq B$ and a non-zero subgroup
 $H\le G$ as in Kemperman's theorem.

If $H$ is infinite or $H=G$ then $(A,B)$ is an elementary pair. Moreover, if
$A$ is not an arithmetic progression, as we can assume, then $(A,B)$ is
elementary of type (III) or (IV); thus, there is a finite subgroup $F\le G$
such that $A$ is contained in a coset of $F$ and $|A+B|\ge|F|-1$. Since
$|A+B|\le |G|-2$, we have $F\neq G$. This shows that $A$ is contained in a
coset of a proper subgroup, and the assertion follows.

Finally, consider the case where $H$ is a finite proper subgroup of $G$. Then
$|A+H|-|A|\le|H|-1$ and we may assume that $A+H=G$; that is, there exist
$g_1\longc g_{n-1}\in G$ such that $A$ is the disjoint union
$A=A_0\cup(g_1+H)\cup\dots\cup(g_{n-1}+H)$. Since $A_0+B_0+H$ has a unique
representation as an element of $\pgh(A)$ and an element of $\pgh(B)$, we
have $B=B_0$. If $A_0$ is not an arithmetic progression, then, as above, the
pair $(A_0,B_0)$ is elementary of type (III) or (IV), and there exists a
non-zero subgroup $F\le H$ such that $A_0$ is contained in a coset of $F$ and
$|A_0+B_0|\ge|F|-1$. The last estimate shows that $F\neq H$: for the
complement of $A+B$ in $G$ is the complement of $A_0+B_0$ in the
corresponding coset of $H$, and if $F=H$ then this complement contains at
most one element. To complete the proof we observe that $A+F\neq G$ and
$|A+F|-|A|\le|F|-1$.
\end{proof}

We now turn to small doubling sets in the elementary abelian groups
$(\Z/2\Z)^r$ and the proof of Theorem \reft{struct}. Our first observation is
that $\mu(A,A)\ge 2$ holds for any subset $A\seq(\Z/2\Z)^r$ with $|A|\ge 2$;
therefore, for $B=A$ Kemperman's condition \refe{kemcon} reduces to
$\pi(2A)=\{0\}$. We need several lemmas.

\begin{lemma}\label{l:antisym}
Let $r\ge 1$ be an integer. If $H\le(\Z/2\Z)^r$ is a subgroup, $h_0\in H$,
and $S\seq H$ is an $h_0$-antisymmetric subset, then $2S=H\stm(h_0+\pi(S))$.
In particular, if $S$ is aperiodic, then $2S=H\stm\{h_0\}$.
\end{lemma}

\begin{proof}
Let $h\in H$. In order for $h\notin 2S$ to hold it is necessary and
sufficient that $(S+h)\cap S=\est$; equivalently, $S+h=H\stm S$, or
$S+h=h_0+S$. The last condition can be re-written as $h+h_0\in\pi(S)$, or
$h\in h_0+\pi(S)$.
\end{proof}

\begin{lemma}\label{l:elemstr}
Let $r\ge 1$ be an integer. For a subset $A\seq(\Z/2\Z)^r$, the pair $(A,A)$
is elementary if and only if $A$ is a shift of an aperiodic antisymmetric
subset of a subgroup of $(\Z/2\Z)^r$; that is, if there exist an element
$g\in(\Z/2\Z)^r$, a non-zero subgroup $H\le(\Z/2\Z)^r$, a subgroup element
$h_0\in H$, and an aperiodic $h_0$-antisymmetric subset $S\seq H$ so that
$A=g+S$.
\end{lemma}

\begin{proof}
It is easily seen that if $H,S,g$, and $h_0$ are as indicated, then $(A,A)$
is elementary of type (I) (if $|A|=1$) or of type (IV) (if $|A|>1$).

On the other hand, if $(A,A)$ is elementary of type (I) then $|A|=1$ and we
can take $g$ to be the element of $A$, set $S=\{0\}$, choose for $H$ an
arbitrary two-element subgroup, and let $h_0$ be the non-zero element of this
subgroup. Next, it is plain that $(A,A)$ cannot be elementary of type (II),
unless $|A|=1$ (there are no non-trivial arithmetic progressions in
$(\Z/2\Z)^r$). Moreover, $(A,A)$ cannot be of type (III) either: for equality
$|g_1+(H_1\cup\{0\})|=|g_2-(H_2\cup\{0\})|$
yields $|H_1|=|H_2|$, whence $|H_1|+|H_2|+1$ is odd. Finally, if $(A,A)$ is
of type (IV), then $A=g_1+H_1=g_2+H_2$ whence $A+g_1=H_1=(g_1+g_2)+H_2$ is a
$(g_1+g_2)$-antisymmetric subset of $H$.
\end{proof}

\begin{lemma}\label{l:smdbel}
Let $r\ge 1$ be an integer and suppose that $A\seq(\Z/2\Z)^r$ satisfies
$\pi(2A)=\{0\}$. Then $|2A|<2|A|$ holds if and only if $(A,A)$ is an
elementary pair.
\end{lemma}

\begin{proof}
We have to show that if $A$ satisfies $|2A|\le 2|A|-1$ and $\pi(2A)=\{0\}$,
then $(A,A)$ is an elementary pair. Without loss of generality, we assume
that $A$ is not contained in a coset of a proper subgroup of $(\Z/2\Z)^r$. If
$(A,A)$ is not elementary then we find $H\lneqq(\Z/2\Z)^r$ as in Theorem
\reft{dual}, as applied with $B=A$. In view of $\oB=\oA$ (we keep using the
notation of Theorem \reft{dual}) and since there is an element in $G/H$ with
a unique representation as a sum of an element of $\oA$ and an element of
$\oB$, we have then $|\oA|=1$; this, however, contradicts the assumption that
$A$ is not contained in a coset of a proper subgroup.
\end{proof}

\begin{proof}[Proof of Theorem \reft{struct}]
For brevity and to avoid double indexation, we write $G=(\Z/2\Z)^r$
throughout the proof.

Suppose that $|2A|<2|A|$ and let $H:=\pi(2A)$. If $|2A|=|H|$, then $A$ is
contained in an $H$-coset and $|A|>|2A|/2=|H|/2$. Assuming now that $|2A|\ge
2|H|$, write $\oA:=\pgh(A)$. Since $|2\oA|<2|\oA|$ (as it follows from
Kneser's theorem) and $\pi(2\oA)=\{0\}$, by Lemmas \refl{smdbel} and
\refl{elemstr} we have $\oA=\oz+\oS$, where $\oz\in G/H$ and $\oS$ is an
aperiodic antisymmetric subset of a non-zero subgroup $\oF\le G/H$. From
$|2A|\ge 2|H|$ we derive that $|\oS|=|\oA|\ge 2$ and in fact $|\oS|>2$, for
any two-element subset of an elementary abelian $2$-group is periodic;
therefore, $|\oF|=2|\oS|\ge 8$. Choose a subgroup $F\le G$ such that
$\pgh^{-1}(\oF)=F\oplus H$ and a subset $S\seq F$ such that $\pgh(S)=\oS$;
thus, $|F|=|\oF|$ and $|S|=|\oS|$. It is easily seen that $S$ is
antisymmetric (if $\oS$ is $\pgh(f)$-antisymmetric for some $f\in F$, then
$S$ is $f$-antisymmetric) and aperiodic (if $S+g=S$ for some $g\in G$, then
$\oS+\pgh(g)=\oS$, hence $\pgh(g)=0,\,g\in H$, and therefore $g=0$ as $g\in
2S\seq F$). Next, from $\oA=\oz+\oS$ it follows that $A\seq z+S+H$, where $z$
is an (arbitrarily selected) element of $G$ such that $\pgh(z)=\oz$.
Furthermore, $2\oA=\oF\stm\{\of\}$ with some $\of\in\oF$ holds by Lemma
\refl{antisym}, whence $2A=(F\oplus H)\stm(f+H)$ for the appropriate choice
of $f\in F$. (Recall that $H=\pi(2A)$.) Finally, in view of
  $$ 2|A| > |2A| = |F||H|-|H| = 2|S||H|-|H| = 2|S+H|-|H| $$
we have $|S+H|-|A|<|H|/2$.

Now suppose that $F,H,S$, and $A$ is as in the second part of the theorem,
and to simplify the notation assume that $A$ is obtained by removing elements
from the \emph{zero} shift of the set $S+H$. For each $s\in S$ we have then
$|(s+H)\stm A|<|H|/2$, whence $|(s+A)\cap H|>|H|/2$ and therefore
$H\seq(s_1+A)+(s_2+A)=s_1+s_2+2A$ for any $s_1,s_2\in S$. It follows that
$s_1+s_2+H\seq 2A$ and hence $2S+H=2A$. Applying Lemma \refl{antisym} we
obtain
  $$ |2A|=(|F|-|\pi(S)|)|H| \le 2|S||H|-|H| = 2|S+H|-|H| < 2|A|, $$
which completes the proof.
\end{proof}

\section*{Appendix: On elementary pairs of types III and IV}

The following two propositions shed some light on the definitions of
elementary pairs of types (III) and (IV). We prove the first proposition
only; the proof of the second is quit similar and we leave it to the reader.

\begin{proposition}\label{p:III}
Suppose that $G$ is an abelian group.
\begin{itemize}
\item[(i)] Let $H$ be a finite subgroup of $G$, and let $(A,B)$ be a
  pair of subsets of $G$ with $\min\{|A|,|B|\}\ge 2$, satisfying \refe{ApB}
  and \refe{kemcon} and such that $A+B$ is a coset of $H$. Then there exist a
  partition $H=H_1\cup H_2\cup\{0\}$ and elements $g_1,g_2\in G$ such that
  $A=g_1+(H_1\cup\{0\})$ and $B=g_2-(H_2\cup\{0\})$.
\item[(ii)] Let $H=H_1\cup H_2\cup\{0\}$ be a partition of a finite
  subgroup $H\le G$ with non-empty subsets $H_1,H_2$, let $g_1,g_2$ be
  elements of $G$, and set $A:=g_1+(H_1\cup\{0\})$ and
  $B:=g_2-(H_2\cup\{0\})$. Then the pair $(A,B)$ satisfies \refe{ApB} and
  \refe{kemcon}, and $A+B$ is a coset of $H$.
\item[(iii)] Let $H,H_1,H_2,g_1,g_2,A,B$ be as in {\rm (ii)} and suppose that
  some group element $c\neq g_1+g_2$ satisfies $\nu_c(A,B)=1$. Write
  $f:=c-g_1-g_2$ (so that $f\in H$), denote by $F$ the subgroup of $H$,
  generated by $f$, and let $\oA:=\pgh[G/F](A)$ and $\oB:=\pgh[G/F](B)$. Then
  $|\oA+\oB|=|\oA|+|\oB|-1$  and $\oA+\oB$ is a coset of $H/F$. Furthermore,
  $A$ is a union of $F$-cosets and an arithmetic progression $A_0$ with
  difference $f$, and $B$ is a union of $F$-cosets and an arithmetic
  progression $B_0$ with difference $f$; moreover, $A_0+B_0$ is an $F$-coset,
  $|A_0+B_0|=|A_0|+|B_0|-1,\,\mu(A_0,B_0)=1$, and $A_0+B_0$, considered as an
  element of $G/F$, has a unique representation as an element of $\oA$ and an
  element of $\oB$.
\end{itemize}
\end{proposition}

\begin{proof}[Proof of Proposition \refp{III}]
(i) We assume without loss of generality that $A,B\seq H$, so that $A+B=H$.
Since $\pi(A+B)=H\neq\{0\}$ in view of $\min\{|A|,|B|\}\ge 2$, we have
$\mu(A,B)=1$ by \refe{kemcon}. Shifting $A$ and $B$, if necessary, we assume
that $0\in A\cap B$ and $\nu_0(A,B)=1$. Write $H_1:=A\stm\{0\}$ and
$H_2:=-B\stm\{0\}$; then $|H_1|+|H_2|=|A|+|B|-2=|A+B|-1=|H|-1$ by \refe{ApB}
and the theorem of Kemperman-Scherk, and $H_1\cap H_2=\est$: for $h\in
H_1\cap H_2$ leads to the ``forbidden'' representation $0=h+(-h)$.

(ii) See a remark in Section \refs{kempth}, following the definition of
elementary pairs.

(iii) We can assume that $g_1=g_2=0$. Fix $a_0\in A$ and $b_0\in B$ with
$c=a_0+b_0$. Then the two sets $A\stm\{a_0\}-c$ and $-B$ are disjoint, and
since $-B=H_2\cup\{0\}$ is the complement in $H$ of $H_1$, we have
$A\stm\{a_0\}-c\seq H_1$. Indeed, from $|A\stm\{a_0\}|=|H_1|$ it follows that
$A\stm\{a_0\}-c=H_1$, which can be rewritten as
$A\stm\{a_0\}=(A+c)\stm\{c\}$. This shows that shifting $A$ by $c$ results in
the same effect as replacing $a_0$ with $c$. Thus, we can write
 $A=A_0\cup A_1\cup\dotsb\cup A_{m-1}$, where $m\ge 1$ is an integer,
$A_i+F=A_i$ for $i\in[1,m-1]$, and $A_0$ is an arithmetic progression with
difference $c$. (Recall, that $F$ is the subgroup, generated by $c$.)
Similarly, we have $B=B_0\cup B_1\cup\dotsb\cup B_{n-1}$, where $n\ge 1$ is
an integer, $B_i+F=B_i$ for $i\in[1,n-1]$, and $B_0$ is an arithmetic
progression with difference $c$. Clearly, the numeration can be so chosen
that if $A_0$ is a full coset of $F$, then $a_0\in A_0$, and similarly if
$B_0$ is a full coset of $F$, then $b_0\in B_0$.

From $|A|+|B|=|H|+1\equiv 1\pmod{|F|}$ we derive that $|A_0|+|B_0|\equiv
1\pmod{|F|}$, so that in fact $|A_0|+|B_0|=|F|+1$, and it follows that
$A_0+B_0$ is an $F$-coset, $|A_0+B_0|=|A_0|+|B_0|-1$, and $\mu(A_0,B_0)=1$.
Since
\begin{multline*}
  |\oA| + |\oB| = (|A|+|F|-|A_0|)/|F| + (|B|+|F|-|B_0|)/|F| \\
     =(|A|+|B|-1)/|F| - (|A_0|+|B_0|-1)/|F| + 2 = |H/F| + 1,
\end{multline*}
we have $\oA+\oB=H/F$ and $|\oA+\oB|=|\oA|+|\oB|-1$.

It remains to show that $A_0+B_0$ has a unique representation as a sum of an
element of $\oA$ and an element of $\oB$. To this end we recall that if $A_0$
is a coset of $F$, then $a_0\in A_0$, and we claim now that this conclusion
stays true if $A_0$ is \emph{not} a coset of $F$; for, in this latter case we
have $|B_0|=|F|+1-|A_0|\ge 2$, so that if $a_0\in A_i$ with $i\in[1,m-1]$,
then $c=a_0+b_0$ cannot satisfy $\nu_c(A,B)=1$. Similarly, $b_0\in B_0$
holds, and therefore an equality of the form $A_0+B_0=A_i+B_j$ with
$i\in[1,m-1],\,j\in[1,n-1]$ would yield a ``forbidden'' representation of
$c=a_0+b_0$.
\end{proof}

\begin{proposition}\label{p:IV}
Suppose that $G$ is an abelian group.
\begin{itemize}
\item[(i)] Let $H$ be a finite subgroup of $G$, and let $(A,B)$ be a
  pair of non-empty subsets of $G$, satisfying \refe{ApB} and \refe{kemcon}
  and such that $A+B$ is a coset of $H$ with one element removed. Then there
  exist a partition $H=H_1\cup H_2$ with $\pi(H_1)=\pi(H_2)=\{0\}$ and
  elements $g_1,g_2\in G$ such that $A=g_1+H_1$ and $B=g_2-H_2$.
\item[(ii)] Let $H=H_1\cup H_2$ be a partition of a finite subgroup $H\le G$
  with non-empty aperiodic subsets $H_1,H_2$, let $g_1,g_2$ be elements of
  $G$, and set $A:=g_1+H_1$ and $B:=g_2-H_2$. Then the pair $(A,B)$ satisfies
  \refe{ApB} and \refe{kemcon}, and $A+B$ is a coset of $H$ with one element
  removed.
\item[(iii)] Let $H,H_1,H_2,g_1,g_2,A,B$ be as in {\rm (ii)} and suppose that
  some group element $c$ satisfies $\nu_c(A,B)=1$. Write $f:=c-g_1-g_2$ (so
  that $f\in H$), denote by $F$ the subgroup of $H$, generated by $f$, and
  let $\oA:=\pgh[G/F](A)$ and $\oB:=\pgh[G/F](B)$. Then
  $|\oA+\oB|=|\oA|+|\oB|-1$ and $\oA+\oB$ is a coset of $H/F$. Furthermore,
  $A$ is a union of $F$-cosets and an arithmetic progression $A_0$ with
  difference $f$, and $B$ is a union of $F$-cosets and an arithmetic
  progression $B_0$ with difference $f$; moreover, $A_0+B_0$ is an $F$-coset
  with one element removed, $|A_0+B_0|=|A_0|+|B_0|-1$, and $A_0+B_0+F$,
  considered as an element of $G/F$, has a unique representation as an
  element of $\oA$ and an element of $\oB$.
\end{itemize}
\end{proposition}

\section*{Acknowledgement}
The author is grateful to the referee for a careful reading of the paper and
a number of valuable remarks.

\end{document}